\documentclass{amsart}
\usepackage{amssymb,euscript,graphicx,psfrag,rotating}
\renewcommand{\hat}{\widehat}
\renewcommand{\phi}{\varphi}

\newcommand{\mathscr}{\EuScript}

\newcommand{\inv}{^{-1}}

\DeclareMathOperator{\tr}{tr}
\DeclareMathOperator{\grad}{grad}

\DeclareMathOperator{\codim}{codim}

\newcommand{\gr}{Grassmannian}
\newcommand{\gmle}{GE}
\newcommand{\gm}{Grassmann manifold}
\newcommand{\sa}{self-$\sig$-adjoint}

\newcommand{\R}{\mathbb{R}} %real numbers
\newcommand{\C}{\mathbb{C}} %complex numbers
\newcommand{\F}{\mathbb{F}} %real or complex numbers
\newcommand{\Gr}{\mathrm{G}} %Grassmannian
\newcommand{\pos}{\mathrm{Pos}}	%positive definite
\newcommand{\sas}{{\Gr(m,r)}} %sample space
\newcommand{\pas}{\pos(m)} %parameterspace
\newcommand{\ls}{{\F^m}} %linear space
\newcommand{\mat}{\F^{m\times m}} %square matrices
\newcommand{\pr}{\mathbb{P}} %projective
\newcommand{\ps}{{\pr^{m-1}}} %projective space
\newcommand{\G}{\mathcal{G}} %Grassmann distribution

\newcommand{\sig}{\Sigma} %typical parameter
\newcommand{\E}{\mathcal{E}} %ellipsoid
\newcommand{\vol}{\mathrm{vol}} %volume
\newcommand{\T}{\mathrm{T}} %tangent space

\newcommand{\ton}{1,\ldots,n}
\newcommand{\listU}{U_1,\ldots,U_n}
\newcommand{\listd}{d_1,\ldots,d_n}

\newcommand{\ineq}{\hbox{\rm(\ref{ineq0}--\ref{ineq2})}}
\newcommand{\Ineq}{\hbox{\rm(\ref{Ineq0}--\ref{Ineq2})}}

\newtheorem{theorem}{Theorem}
\newtheorem{lemma}{Lemma}
\newtheorem{corollary}{Corollary}
\newtheorem{proposition}{Proposition}
\theoremstyle{remark}

\newtheorem{example}{Example}
\begin{document}
\title{Grassmannian Estimation}
\author{Claude Auderset, Christian Mazza and Ernst A. Ruh}
\thanks{Supported by Swiss National Science Foundation Grant
20-05811.99 and 20-67619.02}

\begin{abstract}
This paper discusses the family of distributions on the Grassmannian $\sas$ of the linear span of $r$ central normal vectors in $\R^m$ or $\C^m$, parametrized by the covariance matrix (up to a positive factor). Our main result is an existence and uniqueness criterion for the maximum likelihood estimate of a sample in $\sas$, based on convexity and asymptotic properties of the log-likelihood. By coupling methods of algebraic geometry and linear programming, we show that almost all samples of size $n>m^2/r(m-r)$ in $\sas$ have a unique MLE.

In the real case, a new, unexpected phenomenon takes place for some values $1<r<m$,  which does not occur in the angular Gaussian case $r=1$. Random samples of some critical size in $\sas$ may have a unique estimate or not, with a positive probability in either case.
\end{abstract}
\maketitle
\section{Introduction}
As stated in \cite{rahman}, the current data deluge
inundating science is remarkable for the rapid proliferation
in new data type. Typical examples are directions
in 
$\R^n$ or elements of the Grassmann manifold $\sas$ of all vector subspaces
of dimension $r$ of $\R^m$ ($0<r<m$), as introduced in \cite{chikuse1}.
Being of increasing importance in practical situations
(see e.g.  \cite{dryden}, \cite{mardia}, \cite{4}, \cite{5},
\cite{6}, \cite{rahman} or \cite{8}), there is a strong need for 
studying various classical inference problems, like for
example maximum likelihood estimation. To deal with
these problems, one can in most cases reparametrize
the manifold and recast the inference problem in
some Euclidean space. However, this can have the effect
of hiding intrinsic geometric properties of the
statistical relevant objects (see below).

A typical example is obtained when dealing with
$\sas$ when $r=1$, the set of axes or directions
in $\R^m$, see e.g. \cite{3} and \cite{9}.
In \cite{kent}, the manifold is endowed
with the {\it angular Gaussian distribution}, that
is of the law of the random direction obtained by retaining
 only the axis of a multivariate centered
gaussian random vector in $\R^m$ of 
covariance matrix $\Sigma$.
Kent and Tyler \cite{kent} derived sufficient conditions for the
existence of the maximum likelihood estimator (MLE)
based on an i.i.d. sample by working on
$\R^{m-1}$; the angular Gaussian distribution
is then equivalent to the Cauchy law. The mathematical
analysis can then be performed in $\R^{m-1}$, at the
cost of loosing nice properties of the problem.
In \cite{amr}, the whole picture was obtained
using mainly convexity. The parameter space 
${\rm Pos}(m)$ consists of positive definite self adjoint
matrices of determinant 1, which is considered as
a Riemanian manifold with a natural metric. The results
derived in \cite{amr} make strong use of this manifold
structure, of the particular form of the log likelihood function
and of the geometric link between the parameter space ${\rm Pos}(m)$
and the sample space $\sas$, $r=1$. Interestingly, the estimated scatter matrix
plays a fundamental r\^ole for multivariate nonparametric tests, where it is known as the 
{\em Tyler's transformation matrix}, see e.g. \cite{oja}, or in finance where the maximum likelihood
estimator is used to fit financial data, see \cite{bouchaud}.

When $r$ is arbitrary, we obtain random subspaces
by retaining only the linear span $U=<x_1,\cdots,x_r>$ 
of an i.i.d. sample of $r$ multivariate centered gaussian random
vectors of covariance matrix $\Sigma\in {\rm Pos}(m)$. The
law of this random subspace has been considered previously
in the literature and has been termed as the {\it matrix angular Gaussian distribution} (see e.g. \cite{chikuse1}, \cite{chikuse2} or
\cite{chikuse3}); however, basic questions
like the existence of the MLE remain unexplored. 

We will show that
a new phenomenon emerges: In most statistical settings,
the MLE based on some sample $u_1,\cdots,u_n$ exists
with probability one when the size $n$ is larger that
a critical value $n_c$ and does not exist with
probability one when $n\le n_c$, like for example
in the angular Gaussian case with $r=1$ (see 
e.g. \cite{amr}).
In the Grassmannian  setting,  we show in Example 2 of Section 3 that there are sizes $n$
such that the MLE exists with positive probability
and does not exist with positive probability (see e.g.
\cite{albert} where a similar phenomenon occurs in logistic
regression).

Section 2 introduces the Grassmannian statistical model
and the related likelihood function. Section 3 considers
the problem of existence and uniqueness of the
Grassmannian maximum likelihood estimate (GE). Our main
results, Theorems 1 and 2 give necessary and sufficient conditions
for the existence of a unique GE. The geometrical setting
is illustrated in Examples 1 and 2. Section 4  provides fundamental
properties of the likelihood function like its convexity when
restricted 
to the geodesics of ${\rm Pos}(m)$. This nice property is then
used to prove Theorems 1 and 2. Theorem 4 of Section 5 shows finally that
the GE of almost all samples of size $n$ is unique when
$$n > \frac{m^2}{r(m-r)}.$$

%%%%%%%%%%%%%%%%%%%%%%%%%%%%%%%%%%%%%%%%%%%%%%%%%%%%%%%%%%%%%%%%%%%%%%%%%%%%%%%%%%%%%%%%%%%%%%%%%%%%%%%%%%%%%%%%%%%%%%%%%%%%%%%%%%%%%%%%%%
\section{The \gr\ statistical model}
\subsection{Grassmannian distributions}
We present two versions of the \gr\ model, real or complex. To treat them in parallel, we set $\F=\R$ or $\C$, and denote by $A^*\in\F^{s\times r}$ the adjoint of a matrix $A\in\F^{r\times s}$, i.e., the transpose of $X$ if $\F=\R$ and the complex conjugate of the transpose of $A$ if $\F=C$. A square matrix $\sig\in\mat$ is self-adjoint when $\sig=\sig^*$, i.e., symmetric if $\F=\R$ and Hermitian if $\F=\C$.

Let $x_1,\ldots, x_r\in\ls$ be i.i.d. random vectors in $\ls$ with central normal distribution of positive definite  self-adjoint covariance matrix~$\sig$. The density of the normal law is $\exp(x^*\sig\inv x/2)$ ($x\in\ls$) up to a constant factor in both the real and the complex case. We define the \emph{\gr\ distribution of parameter $\sig$} as the law of the linear span $\langle x_1,\ldots, x_r\rangle$ of these vectors in $\ls$. It is a Borel probability measure $\G_\sig$ on the \emph{\gm\ $\sas$} of all vector subspaces of dimension $r$ of $\ls$ ($0<r<m$). The parameter $\sig$ of a \gr\ distribution $\G_\sig$ is defined up to a a positive factor only. We remove this indeterminacy by requiring the determinant of $\sig$ to be $1$. So, we parametrize the \gr\ distributions by the space $\pas$ of positive definite self-adjoint matrices $\sig\in\mat$ of determinant $1$.

Given a regular matrix $A\in\mat$, the random vectors $Ax_1,\ldots, Ax_r$ are i.i.d.\ with central normal law of covariance matrix $A\sig A^*$. Hence, the image measure of $\G_\sig$ under the transformation of $\sas$ given by
$AU=\{Ax\mid x\in U\}$ for $U\in\sas$ is
\begin{equation}\label{equivariance}
A\G_\sig=\G_{A\sig A^*}.
\end{equation}
In fact, \emph{the \gr\ statistical model $(\G_\sig)_{\sig\in\pas}$ is the unique family of Borel probability measures on $\sas$ indexed by $\pas$ enjoying the  equivariance property~\eqref{equivariance} for all matrices $A\in\mat$ of determinant $1$.} To see this, observe that condition \eqref{equivariance} implies the invariance of $\G_\sig$ under the group of invertible matrices $A$ of determinant~$1$ such that $A\sig A^*=\sig$.  As this group is compact and acts continuously and transitively on $\sas$, there is a unique Borel probability measure  on $\sas$ which is invariant under it, namely $\G_\sig$.

Let us represent a point $U\in\sas$ as the linear span $U=\langle x_1,\ldots,x_r\rangle$ of
linearly independent vectors $x_1,\ldots, x_r$ of $U$ or, equivalently, as the range $U=\langle X\rangle$ of the matrix $X=(x_1,\ldots,x_r)$ of rank $r$. Then, a computation shows that the density, or Radon-Nikodym derivative, of the \gr\ distribution $\G_\sig$ ($\sig\in\pas$) with respect to the uniform distribution $\G_I$ on $\sas$ ($I$ = identity matrix) is given by
\begin{equation}\label{density}
\frac{d\G_\sig}{d\G_I}(\langle
X\rangle)=
\left(\frac{\det(X^*X)}{\det(X^*\sig\inv X)}\right)^{i_\F m/2},
\end{equation}
where $i_\F=\dim_\R(\F)$ (see \cite{chikuse1} for the real case).
The meaning of this formula is perhaps more apparent in the form
\[
\frac{d\G_\sig}{d\G_I}(U)
=
\left(\frac{\vol(\E_I\cap U)}{\vol(\E_\sig\cap U)}\right)^m
\qquad(U\in\sas),
\]
where $\E_\sig=\{x\in\ls\mid x^*\sig\inv x\le1\}$ denotes the ellipsoid
associated to $\sig$ ($\E_I$ = unit ball), and $\vol$ the
Lebesgue measure on $U$.

When $r=1$, the \gr\ distribution $\G_\sig$ is known as the (real or complex) \emph{angular Gaussian distribution} of parameter $\sig\in\pas$ on the projective space $\ps=\Gr(m,1)$ (see \cite{amr}). For any $0<r<m$, the \gm\ $\sas$ can be viewed as the space of projective subspaces of dimension $r-1$ of $\ps$ by identifying a vector $r$-subspace $U$ of $\ls$ with the projective subspace $\{y\in\ps\mid y\subseteq U\}$.
\emph{In this projective interpretation, the \gr\ distribution $\G_\sig$ on $\sas$ is the law of the projective span of i.i.d.\ random points $y_1,\ldots,y_r$ of $\ps$ with angular Gaussian distribution of parameter $\sig$.}

%%%%%%%%%%%%%%%%%%%%%%%%%%%%%%%%%%%%%%%%%%%%%%%%%%%%%%%%%%%%%%%%%%%%%%%%%%%%%%%%%%%%%%%%%%%%%%%%%%%%%%%%%%%%%%%%%%%%%%%%%%%%%%%%%%%%%%%%%%%%%%%%%%%%%%%%%%%%%%

%%%%%%%%%%%%%%%%%%%%%%%%%%%%%%%%%%%%%%%%%%%%%%%%%%%%%%%%%%%%%%%%%%%%%%%%%%%%%%%%%%%%%%%%%%%%%%%%%%%%%%%%%%%%%%%%%%%%%%%%%%%%%%%%%%
\subsection{Grassmannian maximum likelihood estimates}
%%%%%%%%%%%%%%%%%%%%%%%%%%%%%%%%%%%%%%%%%%%%%%%%%%%%%%%%%%%%%%%%%%%%%%%%%%%%%%%%%%%%%%%%%%%%%%%%%%%%%%%%%%%%%%%%%%%%%%%%%%%%%%%%%%%%%

Let $P$ be a Borel probability measure on $\sas$. Typically, we think of $P$ as being the empirical measure $(\delta_{U_1}+\cdots+\delta_{U_n})/n$ of a sample $\listU$ in $\sas$, but other cases are of interest too. A parameter $\sig\in\pas$ is called a \emph{Grassmannian (maximum likelihood) estimate} ---abbreviated \gmle\ in the sequel--- of $P$ if it maximizes the log-likelihood 
$\int_\sas\log(d\G_\sig/d\G_I)\,dP$.
It is called a \gmle\ of a sample $\listU\in\sas$ when $P$ is the empirical measure $(\delta_{U_1}+\cdots+\delta_{U_n})/n$. 

For convenience, we shall rather work with the following negative version of the log-likelihood
\begin{align}
\label{log-likelihood}\ell_P(\sig)&=-\frac1{i_\F m}\int_\sas\log(d\G_\sig/d\G_I)\,dP=\int_\sas\ell_U(\sig)\,dP(U),\\
\noalign{\noindent where the (negative) log-density $\ell_U$ is defined by}
\label{log-density}
\ell_U(\sig)&=-\frac1{i_\F m}\log\frac{d\G_\sig}{d\G_I}(U)
=\frac12\log\frac{\det(X^*\sig\inv X)}{\det(X^*X)}
\qquad(U=\langle X\rangle\in\sas).
\end{align}
With this notation, a \gmle\ of $P$ minimizes $\ell_P$.

%%%%%%%%%%%%%%%%%%%%%%%%%%%%%%%%%%%%%%%%%
\section{Existence and uniqueness of the \gr\ estimate}
%
%%%%%%%%%%%%%%%%%%%%%%%%%%%%%%%%%%%%%%%%%%
%
\begin{theorem}\label{main}
A Borel probability measure $P$ on the real or complex Grassmannian $\sas$ has a unique \gmle\  if and only if
\begin{equation}\label{gcond}
\int_\sas\dim(U\cap V)\,dP(U)<\frac rm\dim(V)
\end{equation}
for all nontrivial linear subspaces $V$ of $\ls$ ($0\ne V\ne\ls$).
\end{theorem}
In the case of an empirical measure $P=(\delta_{U_1}+\cdots+\delta_{U_n})/n$,
\begin{corollary}\label{cormain}
A sample $U_1,\ldots,U_n$ in the real or complex Grassmannian $\sas$  has a unique \gmle\ if and only if
\begin{equation}\label{cond}
\frac1n\sum_{i=1}^n\dim(U_i\cap V)<\frac rm\dim(V)
\end{equation}
for all nontrivial linear subspaces $V$ of $\ls$ ($0\ne V\ne\ls$).
\end{corollary}
The proof of the theorem will be presented in the next section. Let us first consider some special cases.
\begin{example}\label{proj}
When $r=1$, $\sas$ is the projective space $\ps$, and the \gr\ distributions are known as angular Gaussian distributions. In this case, $\dim(U_k\cap V)=1$ or $0$ in Corollary~\ref{cormain} according to whether $U_k\subseteq V$ or not. Hence, \emph{the necessary and sufficient condition for a sample of size $n$ in $\ps$ to have a unique angular Gaussian maximum likelihood estimate is that the number of points of the sample contained in a nontrivial vector subspace $V$ of $\ls$ be less than $n\dim(V)/m$} (see \cite{amr} for a more precise result). 

Now, almost all samples in $\ps$ are in general position, ie.,  any nontrivial vector subspace $V$ of $\ls$ contains at most $\dim(V)$ points of the sample. Thus \emph{almost all samples of size $n>m$ in $\ps$ have a unique angular Gaussian maximum likelihood estimate.} This result goes back to \cite{tyler}.
On the other hand, \emph{no samples of size $n\le m$ in $\ps$ have a unique angular Gaussian maximum likelihood estimate} since any point $U\in\ps$ of a sample is, of course, contained in the one-dimensional subspace $V=U$ of $\ls$, so that the condition for the number of points of the sample contained in $V$ to be less than $n\dim(V)/m$ is not satisfied when $n\le m$.
\end{example}
For a Grassmannian $\sas$ which is not a projective space, the situation is more involved, even in the simplest case $m=4$, $r=2$.
\begin{example}\label{lines}
Let  $\listU$ be a sample in the \gm\ $\Gr(4,2)$, viewed as the space of lines in the projective space $\pr^3$. Suppose that the lines $\listU$ are pairwise skew, i.e., $U_i\cap U_j=0$ for $i\ne j$.  Examining case by case all of the possible  values of $\dim(V)$ and $\dim(U_i\cap V)$ in Corollary \ref{cormain}, we find that the sample has a unique \gmle\ if and only if $n>k$, where $k$ is the maximum number of lines of the sample all of which are met by some line $V\in\Gr(4,2)$. Now, given a line $V$, we can choose any number $n$ of pairwise skew lines $\listU\in\Gr(4,2)$ meeting $V$, so that $n=k$. Hence, there are arbitrary large samples of pairwise skew lines not having a unique \gmle. What is needed is a bound for $k$. 

Recall that the lines meeting each of three pairwise skew lines $U_1$, $U_2$ and $U_3$ form a one-dimensional family $\EuScript{F}_1$ of lines on a quadric surface $Q\subset\pr^3$, whereas the other family $\EuScript{F}_2$ of lines on $Q$ consists of the lines meeting every line of $\EuScript{F}_1$. A point of intersection $x$ of a further line $U_4\in\Gr(4,2)$ with the quadric $Q$ determines a line meeting each of the four lines $U_1$, $U_2$, $U_3$ and $U_4$, namely the line $V\in\EuScript{F}_1$ through~$x$, and vice versa (see Fig.~\ref{fig:quadric}). 
\begin{figure}[thbp]
\psfrag{U1}[]{$U_1$}
\psfrag{U2}[]{$U_2$}
\psfrag{U3}[]{$U_3$}
\psfrag{U4}[]{$U_4$}
\psfrag{V}[]{$V$}
\psfrag{x}[]{$x$}
\centering
\includegraphics{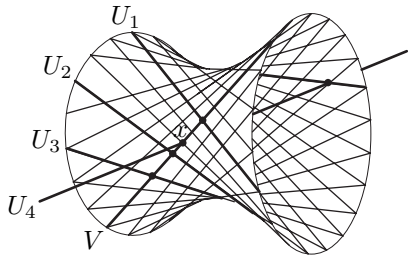}
\caption{Two lines meeting each of four lines}
\label{fig:quadric}
\end{figure}

The number of lines meeting four pairwise skew lines $U_1$, $U_2$, $U_3$ and $U_4$ is thus
\begin{itemize}
\item $2$ if $U_4$ meets $Q$ transversally,
\item $0$ if $U_4$ does not meet $Q$, which can occur only when $\F=\R$,
\item $1$ if $U_4$ is tangent to $Q$,
\item infinite if $U_4$ lies on $Q$,  in which case $U_4\in\EuScript{F}_2$ so that every line meeting $U_1$, $U_2$ and $U_3$ necessarily meets $U_4$ too.
\end{itemize}
In both the real and the complex case, there are at most two lines $V\in\Gr(4,2)$ meeting each of four pairwise skew lines,  except when the four lines belong to the same family of lines on a smooth quadric. So, almost all samples of size $n$ in  $\Gr(4,2)$ consist of pairwise skew lines of which at most four are intersected by a line $V\in\Gr(4,2)$. We conclude from the criterion above that \emph{almost all samples of size $n>4$ in the real or complex \gm\ $\Gr(4,2)$ have a unique \gmle.}

In the complex case, there is a line meeting each of any four pairwise skew lines $U_1$, $U_2$, $U_3$ and $U_4$ since $U_4$ always meets the quadric $Q$. The same holds if some of the four lines meet together or even coincide. Thus, by Corollary~\ref{cormain},  \emph{no samples of size $n\le4$ in the complex \gm\ $\Gr(4,2)$ have a unique \gmle.}

The situation is different in the real case since $U_4$ need not meet the quadric~$Q$. If we choose four lines at random, there may be a line meeting each of them or not, with a positive probability in both cases. Therefore, \emph{the probability that a random sample of size $n=4$ in the real \gm\ $\Gr(4,2)$ has a unique \gmle\ is positive and $<1$.} 

On the other hand, by Corollary~\ref{cormain}, \emph{no samples of size $n<4$ in the real \gm\ $\Gr(4,2)$ have a unique \gmle\ } since any $n<4$ lines are intersected by some line (in fact, by infinitely many lines).
 \end{example}

\section{Likelihood equation\label{LikelihoodEquation}}

We first introduce notions from linear algebra which
are necessary to settle the likelihood equation on
the symmetric space Pos(m). Consider the scalar product
\begin{equation}\label{scalar}
(x|y)_\sig =x^*\sig^{-1}y
\qquad(x,y\in\ls)
\end{equation}
associated to a parameter $\sig\in\pas$. 
We denote by $\pi_U(\sig)$ the \emph{$\sig$-orthogonal projector} onto a vector
subspace~$U$ of $\ls$. It is the linear map
$\pi_U(\sig):\ls\to\ls$ defined by $\pi_U(\sig)u=u$ if
$u\in U$, and $\pi_U(\sig)v=0$ if $v\in\ls$ is \emph{$\sig$-orthogonal} to $U$,
i.e., $(v|u)_\sig=0$ for all $u\in U$. In matrix notation,
\begin{equation}\label{projection}
\pi_U(\sig)=X(X^*\sig\inv X)\inv X^*\sig\inv,
\
\end{equation}
where $U=\langle X\rangle$ is the range of $X\in\F^{m\times r}$. We call a matrix $A\in\mat$ \emph{\sa} if $(Ax|y)_\sig=(x|Ay)_\sig$ for all
$x,y\in\ls$ or, equivalently, if it coincides with its \emph{$\sig$-adjoint} $\sig A^*\sig\inv$.

The parameter space $\pas$ is a Riemannian manifold, in fact a symmetric space. Its tangent space $\T_\sig$ at $\sig\in\pas$ consists of the \sa\ matrices $v\in\mat$ of trace zero, and the Riemannian metric is defined by the scalar products
\begin{equation}\label{metric}
\langle v_1,v_2\rangle=\tr(v_1v_2)
\qquad(v_1,v_2\in\T_\sig)
\end{equation}
on the tangent spaces $\T_\sig$, where $\tr(A)$ denotes the trace of a matrix~$A$. The geodesic $\gamma:\R\to\pas$ of velocity $v\in\T_\sig$ issuing
from $\sig\in\pas$ is
\begin{equation}\label{geodesic}
\gamma(t)=e^{tv}\sig(e^{tv})^*=e^{2tv}\sig
\qquad(t\in\R),
\end{equation}
where $e^{tv}$ denotes the matrix exponential.

Deriving the expression \eqref{log-density} along a geodesic~\eqref{geodesic} and using the matrix form~\eqref{projection} of the $\sig$-orthogonal projector $\pi_U(\sig)$ onto $U$, we find
the gradient (with respect to the Riemannian metric~\eqref{metric} defined above) of  
the log-density
\begin{align}\label{gradient}
\grad\ell_U(\sig)&=\frac r mI-\pi_U(\sig)
\qquad(U\in\sas, \sig\in\pas),\\
\noalign{\noindent and the covariant derivative of $\grad\ell_U$ in the direction of
$v\in\T_\sig$}
\label{covariant}
\nabla_v\grad\ell_U(\sig)&=
\pi_U(\sig)v(I-\pi_U(\sig))+(I-\pi_U(\sig))v\pi_U(\sig).
\end{align}
By integrating these formulas with respect to $P$, and interchanging integration and derivation by means of the Lebesgue dominated convergence theorem,
we get the gradient of the log-likelihood \eqref{log-likelihood}
\begin{align}\label{gradient1}
\grad\ell_P(\sig)&=\frac r mI-\int_\sas\pi_U(\sig)\,dP(U)
\qquad(\sig\in\pas),\\
\noalign{\noindent and its covariant derivative}
\label{covariant1}
\nabla_v\grad\ell_P(\sig)&=\int_\sas
[\pi_U(\sig)v(I-\pi_U(\sig))+(I-\pi_U(\sig))v\pi_U(\sig)]dP(U).
\end{align}
A function $f$ on $\pas$ is called \emph{convex} if its restriction
$f(\gamma(t))$ ($t\in\R$) to any geodesic $\gamma$ is convex in the  
usual sense. This amounts to saying that the Hessian $\nabla^2f$ is positive
semi-definite, i.e.,
$\nabla^2_vf(\sig)=\langle \nabla_v\grad f(\sig),v\rangle\ge0$ for all
$\sig\in\pas$ and $v\in\T_\sig$, since
\begin{equation}\label{secondder}
\frac {d^2}{dt^2}f(\gamma(t))
=(\nabla^2_vf)(\gamma(t))=\langle \nabla_v\grad f(\gamma(t)),v\rangle,
\end{equation}
where $v$ is the velocity of the geodesic $\gamma$.
\begin{proposition}\label{convex}
The log-likelihood function $\ell_P$ is convex. More  
precisely,
its restriction $\ell_P(\gamma(t))$ ($t\in\R$) to a geodesic $\gamma$  
is either strictly convex or affine linear. The latter case occurs if and
only if $v(U)\subseteq U$ for $P$-almost all $U\in\sas$, where 
$v$~is the velocity of the geodesic.
\end{proposition}
\proof

The convexity can be obtained directly by proceeding as in \cite{amr}. On the other hand,
one can use the fact that the log-likelihood function is a {\it Busemann function}
for the symmetric space Pos(m) (see e.g. \cite{fluege}), and convexity follows.
\endproof

As the log-likelihood function $\ell_P$ is convex, its minima are exactly the zeroes  
of its gradient hence, by formula~\eqref{gradient1},
\begin{theorem}
A parameter $\sig\in\pas$ is a \gmle\ of a Borel probability measure~$P$ on $\sas$ if and only if it satisfies the
\emph{maximum likelihood equation}
\begin{equation}\label{LikelihoodFormula}
\int_\sas\pi_U(\sig)\,dP(U)=\frac rm I.
\end{equation}
\end{theorem}

\noindent {\it Proof of Theorem \ref{main}}
\bigskip

\noindent One can either proceed as in \cite{amr}, or use the fact that
the log-likelihood functions is a Busemann function of the 
symmetric space Pos(m), see e.g. \cite{fluege}. The maximum likelihood estimator is then
the barycenter of the related probability measure on the Grassman manifold,
viewed as an orbit in the Tits boundary.Theorem \ref{main} then follows from
Proposition 6.2 of \cite{kapovich}.
\endproof
\section{The linear programming bound}
In order to apply the criteron of Corollary~\ref{cormain} for the existence and uniqueness of the \gmle\ of a sample, we must first answer the following question. 

\emph{Given vector subspaces $U_1,\ldots,U_n\in\sas$ of dimension $r$ of $\ls$ and integers $d_1,\ldots,d_n\ge0$, on what conditions is there a vector subspace $V\in\Gr(m,s)$ of dimension $s$ of $\ls$ such that $\dim(U_k\cap V)=d_k$ for $k=1,\ldots,n$?} A necessary condition, using methods of algebraic geometry, is given by Proposition~\ref{intersection} below. 

In a second step, we look for all possibilities with $0<s=\dim V<m$ and 
\[
\frac1n\sum_{k=1}^n d_k<\frac{rs}m
\]
using methods of linear programming. This leads to the following.

\begin{theorem}\label{mainbound}
Almost all samples of size 
\[
n>\frac{m^2}{r(m-r)}
\]
in the real or complex \gm\ $\sas$ have a unique \gmle.
\end{theorem}
Our main tool is the Schubert calculus on the Grassmannian $\Gr(m,s)$. In general,  the \emph{Schubert variety} (\cite{hodge}, \cite{fulton}) associated to a Young diagram or partition
\[
\lambda=(\lambda_1,\ldots,\lambda_{m-s})\qquad (s\ge\lambda_1\ge\lambda_2\ge\cdots\ge\lambda_{m-s}\ge0)
\]
with at most $s$ rows and $m-s$ columns and a complete flag 
\[
0=F_0\subset F_1\subset\cdots\subset F_m=\ls
\]
of vector subspaces of $\ls$ is defined as
\[
\Omega_\lambda=\{V\in\Gr(m,s)\mid\dim(F_{m-s+i-\lambda_i}\cap V)\ge i, 1\le i\le s\}.
\]
It is an irreducible algebraic subvariety of codimension $|\lambda|=\lambda_1+\cdots+\lambda_{m-s}$ of the Grassmannian $\Gr(m,s)$ of dimension $s(m-s)$. 

In particular, given $U\in\sas$ and an integer $d$ such that 
\[
\max\{0,r+s-m\}\le d_k\le\min\{r,s\},
\]
the set
\[
S_d(U)=\{V\in\Gr(m,s)\mid \dim(U\cap V)\ge d\}
\]
is the Schubert variety $\Omega_\lambda$ associated to the rectangular Young diagram $\lambda=d^k$ with $k=m+d-r-s$ rows and $d$ columns if we choose the flag in such a way that $F_r=U$.
So, 
\[
\codim S_d(U)=d(m+d-r-s).
\]
\begin{proposition}\label{intersection}
The following property holds for almost all samples $(\listU)$ in the real or complex \gm\ $\sas$.
For any vector subspace $V$ of dimension $s$ of $\ls$,
\begin{align}
\label{ineq0}&\max\{0,r+s-m\}\le d_k\le\min\{r,s\}\text{ for $k=\ton$, and}\\
\label{ineq1}&\sum_{k=1}^nd_k(m+d_k-r-s)\le s(m-s),
\end{align}
where $d_k=\dim(U_k\cap V)$.
\end{proposition}
\noindent\emph{Remark.} The conditions~\eqref{ineq0} and~\eqref{ineq1} are necessary for the existence of a vector subspace $V$ such that $d_k=\dim(U_k\cap V)$ for $k=\ton$. But they are not sufficient, as shown by the example $m=6$, $r=3$, $s=3$, $n=2$, $d_1=d_2=2$. In this case, the inequalities~\eqref{ineq0} and~\eqref{ineq1} are satisfied, although
there is in general no $V\in\Gr(6,3)$ meeting $U_1$ and $U_2$ in subspaces of dimension 2.

To get necessary and sufficient conditions, we need the Schubert calculus. But computations in the Schubert calculus (Littlewood-Richardson coefficients) are algorithmically hard \cite{complexity} so we must content ourselves with Proposition~\ref{intersection}.

\proof[Proof of Proposition~\ref{intersection}] The inequalities~\eqref{ineq0} for $d_k=\dim(U_k\cap V)$ follow from the dimension formula
\[
\dim(U_k\cap V)+\dim(U_k+V)=\dim(U_k)+\dim(V).
\]
The proof  of the rest of the proposition uses standard methods of algebraic geometry. 

Let $d_1,\ldots,d_n$ be arbitrary integers satisfying the inequalities~\eqref{ineq0} and consider the algebraic correspondence
\[
C=\{((\listU),V)\in\sas^n\times\Gr(m,s)\mid\dim(U_k\cap V)\ge d_k, k=1,\ldots,n\}.
\]
The range of $C$ is the whole of $\Gr(m,s)$, and its domain $A_{\listd}$ consists of the samples $(\listU)\in\sas^n$ for which there is some $V\in\Gr(m,s)$ with $\dim(U_k\cap V)\ge d_k$ for $k=\ton$.  Let $((\listU),V)$ be a generic point of $C$. Observe that 
\begin{align*}
C\inv(V)&=\{(\listU)\in\sas^n\mid((\listU),V)\in C\}\\
&=S_{d_1}(V)\times\cdots\times S_{d_n}(U_n),
\end{align*}
where $S_{d_k}(V)=\{U\in\sas\mid \dim(U\cap V)\ge d_k\}$ is a Schubert variety with
\[
\codim S_{d_k}(V)=d_k(m+d_k-r-s)
\]
as explained above for $S_d(U)$.

According to the principle of counting constants \cite{hodge},
\[
\dim A_{\listd}+\dim C(\listU)=\dim \Gr(m,s)+\dim C\inv(V),
\]
where $C(\listU)$ consists of all $V\in\Gr(m,s)$ such that $\dim(U_k\cap V)\ge d_k$ for $k=\ton$, hence
\begin{align*}
\dim A_{\listd}&\le \dim \Gr(m,s)+\dim C\inv(V)\\
&=s(m-s)+\sum_{k=1}^n\bigl(r(m-r)-d_k(m+d_k-r-s)\bigr)\\
&=\dim \sas^n+s(m-s)-\sum_{k=1}^nd_k(m+d_k-r-s).
\end{align*}
This shows that $A_{\listd}$ is a proper algebraic subset of $\sas^n$ if the inequality~\eqref{ineq1} is not satisfied. Let $N_s$ be the union of $A_{\listd}$ where $(\listd)$ runs over all those lists of integers  satisfying the inequalities~\eqref{ineq0} but not the inequality~\eqref{ineq1}, and let $N$ be the union of $N_s$ for $s=0,\ldots,m$. As a finite union of proper algebraic subsets, $N$ is also a proper algebraic subset by the irreducibility of $\sas^n$, hence negligible. Now, take a sample $(\listU)\in\sas^n$ not belonging to~$N$, and any vector subspace~$V$ of any dimension $s$ of $\F^m$. Set $d_k=\dim(U_k\cap V)$ for $k=\ton$, so that $(\listU)\in A_{\listd}$. Then $(\listd)$ must satisfy the inequality~\eqref{ineq1}, otherwise $(\listU)$ would belong to $N$ by the very definition of~$N$. This proves the Proposition.
\endproof

Consider next the set $B(m,r,s)$ of those positive integers $n$ for which there are integers $\listd$ satisfying the inequalities
\begin{align}
\addtocounter{equation}{-2}
&\max\{0,r+s-m\}\le d_k\le\min\{r,s\}\quad\text{for $k=\ton$,}\\
&\sum_{k=1}^nd_k(m+d_k-r-s)\le s(m-s),\\
\label{ineq2}&m\sum_{k=1}^nd_k\ge nrs,
\end{align}
and set $B(m,r)=\bigcup_{s=1}^{m-1}B(m,r,s)$. 
\begin{lemma}\label{lemma1}
Almost all samples of size $n\notin B(m,r)$ in the real or complex \gm\ $\sas$ have a unique \gmle. \end{lemma}
\proof
Suppose that $n\notin B(m,r)$.  According to Proposition~\ref{intersection}, the following holds for almost all $(\listU)\in\sas^n$. For any proper vector subspace $V$ of dimension $s$ of $\ls$, the integers $d_k=\dim(U_k\cap V)$ satisfy the inequalities~\eqref{ineq0} and~\eqref{ineq1}. But they do not satisfy the inequality~\eqref{ineq2} since $n\notin B(m,r)$ hence $n\notin B(m,r,s)$. Thus $m\sum_{k=1}^nd_k<nrs$, which is precisely the condition~\eqref{cond} of Corollary~\ref{cormain} for the sample $\listU$ to have a unique \gmle.
\endproof
\begin{lemma}\label{bound}
For any integers $m,r$ with $0<r<m$, the set $B(m,r)$ is bounded above by $m^2/r(m-r)$.
\end{lemma}
\proof
As $B(m,r)=\bigcup_{s=1}^{m-1}B(m,r,s)$, we first look for an upper bound of $B(m,r,s)$. To this end, we replace the unknowns $\listd$ in the the definition of $B(m,r,s)$ by the number 
\[
n_i=\#\{k\in\{1,\ldots,n\}\mid d_k=i\}
\]
of occurences among $\listd$ of each integer $i$ between $i_0$ and $i_1$, where
\[
i_0=\max\{0,r+s-m\}\quad\text{and}\quad i_1=\min\{r,s\}.
\]
With these new unknowns $n_{i_0},\ldots,n_{i_1}$, the inequations \ineq\ translate into the system of linear inequations
\begin{align}
\label{Ineq0}&n_i\ge0\quad\text{for $i_0\le i\le i_1$},\\
\label{Ineq1}&\sum_{i=i_0}^{i_1}i(m+i-r-s)n_i\le s(m-s),\\
\label{Ineq2}&\sum_{i=i_0}^{i_1}(rs-mi)n_i\le0,
\end{align}
with $n=\sum_{i=i_0}^{i_1}n_i$. So, $B(m,r,s)$ consists of those integers $n$ which decompose into a sum $n=\sum_{i=i_0}^{i_1}n_i$ of integers $n_i$ satisfying the inequalities \Ineq. The maximum of $B(m,r,s)$ (if any) is the solution of the integer linear program
\[
\text{maximize }\sum_{i=i_0}^{i_1}n_i\text{ subject to the constraints~\Ineq.}
\]
Relaxing the integrality condition on $n_i$ yields a usual linear program with real $n_{i_0},\ldots,n_{i_1}$, whose solution is an upper bound of $B(m,r,s)$. Standard methods of linear programming \cite{schrijver} show that the constraints \ineq\ define a bounded polytope whose vertices are of one of the following two types.
\begin{itemize}
\item $n_i=\dfrac{s(m-s)}{i(m+i-r-s)}$ for some $i$ and $n_k=0$ for $k\ne i$.
\item $n_i$ and $n_j$ are the solutions of the system of equations
\[
\left\{
\begin{aligned}
&i(m+i-r-s)n_i+j(m+j-r-s)n_j=s(m-s),\\
&(rs-mi)n_i+(rs-m_j)n_j=0.
\end{aligned}
\right.
\]
and $n_k=0$ for $k\ne i, j$. 
\end{itemize}
Now, routine computations show that the sum $n=\sum_{i=i_0}^{i_1}n_i$ reaches its maximum on vertices of the the first type when $i=\lceil rs/m\rceil$, and on vertices of the second type when $i+1=j=\lceil rs/m\rceil$. It can then be checked that these maxima are bounded above by the quantity $m^2/r(m-r)$.
\endproof

Theorem~\ref{mainbound} immediately follows from Lemma~\ref{lemma1} and~\ref{bound}.

\section{Numerical algorithms}
Let $P$ be a probability measure admitting a unique 
maximum likelihood estimator.
We propose here two algorithms to locate this estimator, using
the geometry of the problem (see Section \ref{LikelihoodEquation}).
 The first one
is a gradient-descent dynamics. The second one
is a faster method which avoids the time consuming steps of the first one.

We look for the solution $\hat\sig$ to the equation  (\ref{LikelihoodFormula}). The Exponential
map ${\rm Exp}_\sig$ from $\T_\sig$ to $\pas$ is given explicitely by ${\rm Exp}_\sig (v) =e^{v}\sig$. Given some $\sig_k$ and $\sig'={\rm Exp}_{\sig_k}(v)$, the idea
is to approximate the gradient $\grad\ell_P(\sig')$ using the parallel transport of 
$\grad\ell_P(\sig_k)+\nabla_v(\sig_k)$. One then
computes the solution $v_{k+1}\in\T_{\sig_k}$ to the linear system
\begin{equation}\label{LoopEquation}
\grad\ell_P(\sig_k)+\nabla_{v_{k+1}}\grad\ell_P(\sig_k)=0.
\end{equation}
The loop is closed
by setting $\sig_{k+1}={\rm Exp}_{\sig_k}(v_{k+1})$.

The step which consists in solving
(\ref{LoopEquation}) is time consuming, so that we propose a faster dynamics:
Given $\sig_k$, we use the geodesic $\gamma_k(t)=e^{2t\grad\ell_P(\sig_k)}\sig_k$, and set
\begin{equation}\label{Numerical}
\sig_{k+1}=\gamma_k(1)=e^{2\grad\ell_P(\sig_k)}\sig_k.
\end{equation}

Our simulations indicate that the sequence $(\sig_k)_{k\ge 0}$ converges toward the
maximum likelihood estimator $\hat\sig^n$. We have performed a simulation study using
$n=50,\ 500,\ 5000$ i.i.d. random samples $\langle X^1\rangle,\cdots,\langle X^n\rangle$, $\langle X^i\rangle\in\Gr(4,2)$, distributed
according to the Grassmannian distribution of parameter $\sig_0$ given by

\begin{center}
\begin{tabular}{l c c c}
	&	&	&\\
1.23943	& 0.53234	& 0.21763	& 0.33038\\
0.53234	& 1.12502	& 0.76236	& 0.20842\\
0.21763	& 0.7626	& 1.52821	& 0.82655\\
0.33038	& 0.20842	& 0.82655	& 1.52298\\
\end{tabular}
\end{center}

The probability measure $P$ is then the empirical distribution
on $\Gr(4,2)$
associated with the random sample.
Our simulations indicate that the maximum likelihood is consistent.
For a random sample of size $n=500$, we found that
the difference between $\sig_0$ and the estimate $\hat\sig^n$
is given by

\begin{center}
\begin{tabular}{l c c c}
	&	&	&\\
0.0282495	& 0.0095817	& 0.0791341	& -0.0819841\\
	& 0.0269432	& 0.1031291	& -0.0447055\\
	&	& 0.1444463	& -0.0051798\\
	& 	& 	& -0.1134552\\
\end{tabular}
\end{center}

For $n=5000$, this difference was given by
\begin{center}
\begin{tabular}{l c c c}
	&	&	&\\
0.01223629	& -0.0100086	& 0.0110916	& -0.0221974\\
	& -0.0209799	& -0.0366614	& -0.0114825\\
	&	& -0.0571491	& 0.0010570\\
	& 	& 	& 0.0380609\\
\end{tabular}
\end{center}

\end{document}